\documentclass[bib]{ba}

\usepackage{bayes}
\usepackage{amsthm}
\usepackage{amsmath}
\usepackage{amssymb}
\usepackage{amsfonts}
\usepackage{mathabx}
\usepackage{bbm}
\usepackage{multirow}
\usepackage[breaklinks,plainpages=false,colorlinks=true,linkcolor=blue,%
  anchorcolor=green,citecolor=blue,filecolor=blue,urlcolor=blue]{hyperref}
\usepackage{graphicx}
\DeclareGraphicsExtensions{.pdf}

\newcommand{\Prob}{\operatorname{P}}
\newcommand{\given}{\operatorname{|}}

\newcommand{\tr}{\operatorname{tr}}

\newcommand{\E}{\mathsf{E}}
\newcommand{\VAR}{\mathsf{VAR}}
\newcommand{\COV}{\mathsf{COV}}
\newcommand{\COR}{\mathsf{COR}}

\newcommand{\argmin}{\operatornamewithlimits{argmin}}

\newcommand{\prior}{\Prob(\mathcal{G})}
\newcommand{\posterior}{\Prob(\mathcal{G}\given\mathcal{D})}
\newcommand{\eprior}{\Prob(\mathcal{G}(\mathcal{E}))}
\newcommand{\eposterior}{\Prob(\mathcal{G}(\mathcal{E})\given\mathcal{D})}

\newtheorem{theorem}{Theorem}[section]
\newtheorem{lemma}{Lemma}[section]

\newtheorem{conjecture}{Conjecture}[section]
\newtheorem{example}{Example}[section]

\begin{document}

\inserttype[]{article}
\author{Marco Scutari}{%
  Marco Scutari\\
  Genetics Institute, University College London, United Kingdom\\
  m.scutari@ucl.ac.uk
}
\title[Priors and Posteriors in Graphical Modelling]
  {On the Prior and Posterior Distributions Used in Graphical Modelling}
\maketitle

\begin{abstract}
  Graphical model learning and inference are often performed using Bayesian
  techniques. In particular, learning is usually performed in two separate
  steps. First, the graph structure is learned from the data; then the
  parameters of the model are estimated conditional on that graph structure.
  While the probability distributions involved in this second step have been
  studied in depth, the ones used in the first step have not been explored
  in as much detail.

  In this paper, we will study the prior and posterior distributions defined
  over the space of the graph structures for the purpose of learning the
  structure of a graphical model. In particular, we will provide a
  characterisation of the behaviour of those distributions as a function of
  the possible edges of the graph. We will then use the properties resulting
  from this characterisation to define measures of structural variability 
  for both Bayesian and Markov networks, and we will point out some of their
  possible applications.


\end{abstract}

Graphical models \citep{pearl,lauritzen} stand out among other classes of
statistical models because of their use of graph structures in modelling and
performing inference on multivariate, high-dimensional data. The close
relationship between their probabilistic properties and the topology of the
underlying graphs represents one of their key features, as it allows an
intuitive understanding of otherwise complex models.

In a Bayesian setting, this duality leads naturally to split model estimation
(which is usually called \textit{learning}) in two separate steps \citep{cowell}.
In the first step, called \textit{structure learning}, the graph structure
$\mathcal{G}$ of the model is estimated from the data. The presence (absence)
of a particular edge between two nodes in $\mathcal{G}$ implies the conditional
(in)dependence of the variables corresponding to such nodes. In the second step,
called \textit{parameter learning}, the parameters $\Theta$ of the distribution
assumed for the data are estimated conditional to the graph structure obtained
in the first step. If we denote a graphical model with $\mathcal{M}$, so that
$\mathcal{M} = (\mathcal{G}, \Theta)$, then we can write graphical model
estimation from a data set $\mathcal{D}$ as
\begin{equation*}
  \Prob(\mathcal{M} \given \mathcal{D}) = \posterior \Prob(\Theta \given \mathcal{G}, \mathcal{D}).
\end{equation*}
Furthermore, following \citet{heckerman}, we can rewrite structure learning as
\begin{equation}
\label{eqn:structlearn}
  \posterior \propto \prior\Prob(\mathcal{D} \given \mathcal{G}).
\end{equation}
The prior distribution $\prior$ and the corresponding posterior distribution
$\posterior$ are defined over the space of the possible graph structures, say
$\mathbf{G}$. Since the dimension of $\mathbf{G}$ grows super-exponentially
with the number of nodes in the graph \citep{harary}, it is common practice
to choose
\begin{align}
\label{eqn:flatprior}
  &\prior = \frac{1}{|\mathbf{G}|}& &\text{for every $\mathcal{G} \in \mathbf{G}$}
\end{align}
as a non-informative prior, and then to search for the graph structure
$\mathcal{G}$ that maximises $\posterior$. Unlike such a \textit{maximum
a posteriori} (MAP) approach, a full Bayesian analysis is computationally
unfeasible in most real-world settings \citep{friedman,koller}. Therefore,
inference on most aspects of $\prior$ and $\posterior$ is severely limited
by the nature of the graph space.

In this paper, we approach the analysis of those probability distributions from
a different angle. We start from the consideration that, in a graphical model,
the presence of particular edges and their layout are the most interesting features
of the graph structure. Therefore, investigating $\prior$ and $\posterior$
through the probability distribution they induce over the set $\mathcal{E}$ of
their possible edges (identified by the set of unordered pairs of nodes in
$\mathcal{G}$) provides a better basis from which to develop Bayesian inference
on $\mathcal{G}$. This can be achieved by modelling $\mathcal{E}$ as a
multivariate discrete distribution encoding the joint state of the edges. Then,
as far as inference on $\mathcal{G}$ is concerned, we may rewrite Equation
\ref{eqn:structlearn} as
\begin{equation*}
  \eposterior \propto \eprior\Prob(\mathcal{D} \given \mathcal{G}(\mathcal{E})).
\end{equation*}
As a side effect, this shift in focus reduces the effective dimension of the
sample space under consideration from super-exponential (the dimension of 
$\mathbf{G}$) to polynomial (the dimension $\mathcal{E}$) in the number of nodes.
The dimension of the parameter space for many inferential tasks, such as the
variability measures studied in this paper, is likewise reduced.

The content of the paper is organised as follows. Basic definitions and notations
are introduced in Section \ref{sec:definitions}. The multivariate distributions
used to model $\mathcal{E}$ are described in Section \ref{sec:distributions}.
Some properties of the prior and posterior distributions on the graph space,
$\eprior$ and $\eposterior$, are derived in Section \ref{sec:properties}. We will
focus mainly on those properties related with the first and second order moments
of the distribution of $\mathcal{E}$, and we will use them to characterise several
measures of structural variability in Section \ref{sec:variability}. These measures
may be useful for several inferential tasks for both Bayesian and Markov networks;
some will be sketched in Section \ref{sec:variability}. Conclusions are summarised
in Section \ref{sec:conclusion}, and proofs for the theorems in Sections 
\ref{sec:distributions} to \ref{sec:variability} are reported in Appendix
\ref{app:proofs}. Appendix \ref{app:numbers} lists the exact values for some
quantities of interest for $\eprior$, computed for several graph sizes.

\section{Definitions and notations}
\label{sec:definitions}

Graphical models \citep{lauritzen,pearl} are a class of statistical models which
combine the rigour of a probabilistic approach with the intuitive representation
of relationships given by graphs. They are composed by a set $\mathbf{X} = \{X_1,
\ldots, X_n\}$ of \textit{random variables} describing the data $\mathcal{D}$
and a \textit{graph} $\mathcal{G} = (\mathbf{V}, E)$ in which each \textit{vertex}
or \textit{node} $v \in \mathbf{V}$ is associated with one of the random variables
in $\mathbf{X}$. Nodes and the corresponding variables are usually referred to
interchangeably. The \textit{edges} $e \in E$ are used to express the dependence
relationships among the variables in $\mathbf{X}$. Different classes of graphs
express these relationships with different semantics, having in common the
principle that graphical separation of two vertices implies the conditional
independence of the corresponding random variables \citep{pearl}. The two examples
most commonly found in literature are \textit{Markov networks} \citep{whittaker,
edwards}, which use undirected graphs \citep[UGs, see][]{diestel}, and
\textit{Bayesian networks} \citep{neapolitan,korb}, which use directed acyclic
graphs \citep[DAGs, see][]{digraphs}. In the context of Bayesian networks, edges
are often called \textit{arcs} and denoted with $a \in A$; we will adopt this
notation as well.

The structure of $\mathcal{G}$ (that is, the pattern of the nodes and the edges)
determines the probabilistic properties of a graphical model. The most important,
and the most used, is the factorisation of the \textit{global distribution} (the
joint distribution of $\mathbf{X}$) into a set of lower-dimensional \textit{local
distributions}. In Markov networks, local distributions are associated with
\textit{cliques} (maximal subsets of nodes in which each element is adjacent to
all the others); in Bayesian networks, each local distribution is associated with
one node conditional on its \textit{parents} (nodes linked by an incoming arc).
In Markov networks the factorisation is unique; different graph structures
correspond to different probability distributions. This is not so in Bayesian
networks, where DAGs can be grouped into \textit{equivalence classes} which are
statistically indistinguishable. Each such class is uniquely identified by the
underlying UG (i.e. in which arc directions are disregarded, also known as
\textit{skeleton}) and by the set of \textit{v-structures} (i.e. converging
connections of the form $v_i \rightarrow v_j \leftarrow v_k$, $i \neq j \neq k$,
in which $v_i$ and $v_k$ are not connected by an arc) common to all elements of
the class.

As for the global and the local distributions, there are many possible choices
depending on the nature of the data and the aims of the analysis. However,
literature have focused mostly on two cases: the \textit{discrete case}
\citep{whittaker,heckerman}, in which both the global and the local distributions
are multinomial random variables, and the \textit{continuous case}
\citep{whittaker,heckerman3}, in which the global distribution is multivariate
normal and the local distributions are univariate (in Bayesian networks) or
multivariate (in Markov networks) normal random variables. In the former, the
parameters of interest $\Theta$ are the \textit{conditional probabilities}
associated  with each variable, usually represented as conditional probability
tables. In the latter, the parameters of interest $\Theta$ are the \textit{partial
correlation coefficients} between each variable and its neighbours in $\mathcal{G}$.
Conjugate distributions (Dirichlet and Wishart, respectively) are then used
for learning and inference in a Bayesian setting.

\section{Multivariate discrete distributions}
\label{sec:distributions}

The choice of an appropriate probability distribution for the set $\mathcal{E}$
of the possible edges is crucial to make the derivation and the interpretation
of the properties of $\mathcal{E}$ and $\mathcal{G}(\mathcal{E})$ easier. We
will first note that a graph is uniquely identified by its edge set $E$ (or by
its arc set $A$ for a DAG), and that each edge $e_{ij}$ or arc $a_{ij}$ is
uniquely identified by the nodes $v_i$ and $v_j$, $i \neq j$ it is incident on.
Therefore, if we model $\mathcal{E}$ with a random variable we have that any
edge set $E$ (or arc set $A$) is just an element of its sample space; and since
there is a one-to-one correspondence between graphs and edge sets, probabilistic
properties and inferential results derived for traditional graph-centric
approaches can easily be adapted to this new edge-centric approach and vice
versa. In addition, if we denote $\mathcal{E} = \{ (v_i, v_j), i \neq j\}$, we
can clearly see that $|\mathcal{E}| = \mathcal{O}(|\mathbf{V}|^2)$. On the other
hand, $|\mathbf{G}| = \mathcal{O}(2^{|\mathbf{V}|^2})$ for UGs and even larger
for DAGs \citep{robinson,harary} and their equivalence classes \citep{perlman}.

We will also note that an edge or an arc has only few possible states:
\begin{itemize}
  \item an edge can be either present ($e_{ij} \in E$) or missing from an
    UG ($e_{ij} \notin E$);
  \item in a DAG, an arc can be present in one of its two possible directions
    ($\overleftarrow{a_{ij}} \in A$ or $\overrightarrow{a_{ij}} \in A$) or
    missing from the graph ($\overleftarrow{a_{ij}} \notin A$ and
    $\overrightarrow{a_{ij}} \notin A$).
\end{itemize}
This leads naturally to the choice of a Bernoulli random variable for the former,
\begin{equation}
\label{eqn:bidef}
  e_{ij} \sim E_{ij} = \left\{
    \begin{aligned}
     1&  &e_{ij} \in E &\text{ with probability $p_{ij}$} \\
     0&  &e_{ij} \notin E &\text{ with probability $1 - p_{ij}$} 
    \end{aligned}
  \right.,
\end{equation}
and to the choice of a Trinomial random variable for the latter,
\begin{equation}
\label{eqn:tridef}
  a_{ij} \sim A_{ij} = \left\{
    \begin{aligned}
     -1&  &\overleftarrow{a_{ij}} \in A &\text{ with probability $\overleftarrow{p_{ij}}$} \\
     0&   &\overleftarrow{a_{ij}}, \overrightarrow{a_{ij}} \notin A
          &\text{ with probability $\mathring{p_{ij}}$}  \\
     1&   &\overrightarrow{a_{ij}} \in A &\text{ with probability $\overrightarrow{p_{ij}}$} 
    \end{aligned}
  \right.,
\end{equation}
where $\overrightarrow{a_{ij}}$ is the arc $v_i \rightarrow v_j$ and
$\overleftarrow{a_{ij}}$ is the arc $v_j \rightarrow v_i$. Therefore, a
graph structure can be modelled through its edge or arc set as follows:
\begin{itemize}
  \item UGs, such as Markov networks or the skeleton and the moral graph of
    Bayesian networks \citep{pearl}, can be modelled by a \textit{multivariate
    Bernoulli random variable};
  \item directed graphs, such as the DAGs used in Bayesian networks, can be
    modelled by a \textit{multivariate Trinomial random variable}.
\end{itemize}

In addition to being the natural choice for the respective classes of graphs,
these distributions integrate smoothly with and extend other approaches present
in literature. For example, the probabilities associated with each edge or arc
correspond to the \textit{confidence coefficients} from \citet{friedman} and the
\textit{arc strengths} from \citet{imoto}. In a frequentist setting, they have
been estimated using bootstrap resampling \citep{efron}; in a Bayesian setting,
Markov chain Monte Carlo (MCMC) approaches \citep{friedman2,melancon} have
been used instead.

\subsection{Multivariate Bernoulli}
\label{sec:mvber}

Let $B_1, \ldots, B_k$, $k \in \mathbb{N}$ be Bernoulli random variables with
marginal probabilities of success $p_1, \ldots, p_k$, that is $B_i \sim Ber(p_i)$,
$i = 1, \ldots, k$. Then the distribution of the random vector $\mathbf{B} =
[B_1, \ldots, B_k]^T$ over the joint probability space of $B_1, \ldots, B_k$ is
a \textit{multivariate Bernoulli random variable} \citep{krummenauer}, denoted
as $Ber_k(\mathbf{p})$. Its probability function is uniquely identified by the
parameter collection
\begin{equation*}
  \mathbf{p} = \left\{ p_I : I \subseteq \{1, \ldots, k \},\, I \neq \varnothing \right\},
\end{equation*}
which represents the \textit{dependence structure} among the $B_i$ in terms of
simultaneous successes for every non-empty subset $I$ of elements of $\mathbf{B}$.
Other characterisations and fundamental properties of the multivariate Bernoulli
distribution can be found in \citet{kotz}. \citet{bivdist} focus on the bivariate
models specific to $Ber_2(\mathbf{p})$. Additional characterisations and results
specific to particular applications can be found in \citet[][variable
selection]{george}, \citet[][longitudinal studies]{farrell}, \citet[][combinatorial
optimisation]{rubinstein} and \citet[][clinical trials]{drugs}, among others.

From literature we know that the expectation and the covariance matrix of $\mathbf{B}$
are immediate extensions of the corresponding univariate Bernoulli ones;
\begin{align*}
  &\E(\mathbf{B}) = [ p_1, \ldots p_k]^T& &\text{and}&
  &\COV(\mathbf{B}) = [ \sigma_{ij} ] = p_{ij} - p_{i}p_{j}.
\end{align*}
In particular, the covariance matrix $\Sigma = [\sigma_{ij}]$ has some interesting
numerical properties. From basic probability theory, we know its diagonal elements
$\sigma_{ii}$ are bounded in the interval $\left[0, \frac{1}{4}\right]$; the
maximum is attained for $p_i = \frac{1}{2}$, and the minimum for both $p_i = 0$ and
$p_i = 1$. For the Cauchy-Schwarz theorem then $|\sigma_{ij}| \in \left[0, \frac{1}{4}\right]$.
As a result, we can derive similar bounds for the eigenvalues $\lambda_1, \ldots,
\lambda_k$ of $\Sigma$, as shown in the following theorem.
\begin{lemma}
\label{thm:mvebereigen}
  Let $\mathbf{B} \sim Ber_k(\mathbf{p})$, and let $\Sigma$ be its covariance
  matrix. Let $\lambda_i$, $i = 1, \ldots, k$ be the eigenvalues of $\Sigma$.
  Then
  \begin{align*}
    &0 \leqslant \sum_{i=1}^k \lambda_i \leqslant \frac{k}{4}&
    &\text{and}&
    &0 \leqslant \lambda_i \leqslant \frac{k}{4}.
  \end{align*}
\end{lemma}
\begin{proof}
See Appendix \ref{app:proofs}.
\end{proof}

These bounds define a closed convex set in $\mathbb{R}^k$, described by the
family
\begin{equation*}
  \mathcal{L} = \left\{ \Delta^{k-1}(c) : c \in \left[ 0, \frac{k}{4} \right]\right\}
\end{equation*}
where $\Delta^{k-1}(c)$ is the non-standard $k-1$ simplex
\begin{equation}
\label{eq:simplex}
  \Delta^{k-1}(c) = \left\{ (\lambda_1, \ldots, \lambda_k) \in \mathbb{R}^k :
  \sum_{i=1}^k \lambda_i = c, \lambda_i \geqslant 0\right\}.
\end{equation}

\subsection{Multivariate Trinomial}
\label{sec:mvtri}

Construction and properties of the multivariate Trinomial random variable are
similar to the ones illustrated in the previous section for the multivariate
Bernoulli. For this reason, and because it is a particular case of the
multivariate multinomial distribution, the multivariate Trinomial distribution
is rarely the focus of research efforts in literature. Some of its fundamental
properties are covered either in \citet{kotz} or in monographs on contingency
tables analysis such as \citet{fienberg}. 

Let $T_1, \ldots, T_k$, $k \in \mathbb{N}$ be Trinomial random variables assuming
values $\{-1, 0, 1\}$ and denoted as $T_i \sim Tri\left(p_{i(-1)}, p_{i(0)},
p_{i(1)}\right)$ with $p_{i(-1)} + p_{i(0)} + p_{i(1)} = 1$. Then the distribution
of the random vector $\mathbf{T} = [T_1, \ldots, T_k]^T$ over the joint probability
space of $T_1, \ldots, T_k$ is a \textit{multivariate Trinomial random
variable}, denoted as $Tri_k(\mathbf{p})$. The parameter collection $\mathbf{p}$
which uniquely identifies the distribution is
\begin{equation*}
  \mathbf{p} = \left\{ p_{I(T)} : I \subseteq \{1, \ldots, k\},\,
    T \in \bigtimes_{i = 1}^{|I|} \{-1, 0, 1\},\,
    I \neq \varnothing \right\}
\end{equation*}
and the reduced parameter collection we will need to study its first and second
order moments is
\begin{equation*}
  \mathbf{\tilde{p}} = \left\{ p_{ij(T)} : i,j = 1, \ldots, k,\,
    T \in \{-1, 0,1\}^2 \right\}.
\end{equation*}

From the definition, we can easily derive the expected value and the variance
of $T_i$,
\begin{align*}
  \E(T_i) &= p_{i(1)} - p_{i(-1)} \\
  \VAR(T_i) &= p_{i(1)} + p_{i(-1)} - \left[\, p_{i(1)} - p_{i(-1)}\, \right]^2 
\end{align*}
and the covariance between two variables $T_i$ and $T_j$,
\begin{align*}
  \COV(T_i, T_j) &= \left[\, p_{ij(1,1)} - p_{i(1)}p_{j(1)} \,\right] + 
                    \left[\, p_{ij(-1,-1)} - p_{i(-1)}p_{j(-1)}  \,\right] -  \notag \\
                 &\qquad    - \left[\, p_{ij(-1,1)} - p_{i(-1)}p_{j(1)}  \,\right] 
                     - \left[\, p_{ij(1,-1)} - p_{i(1)}p_{j(-1)}  \,\right].
\end{align*}

Again, the diagonal elements of the covariance matrix $\Sigma$ are bounded. This
can be proved either by solving the constrained maximisation problem
\begin{align*}
  &\max_{p_{i(1)}, p_{i(-1)}} \VAR(T_i)&  
  &\text{s.t.}& &p_{i(1)} \geqslant 0,\; p_{i(-1)} \geqslant 0,\; p_{i(1)} + p_{i(-1)} \leqslant 1
\end{align*}
or as an application of the following theorem by \citet{discvar}.
\begin{theorem}
  If a discrete random variable $X$ can take values only in the segment $[x_1, x_n]$
  of the real axis, the maximum standard deviation of $X$ equals $\frac{1}{2}(x_n - x_1)$.
  The maximum is reached if $X$ takes the values $x_1$ and $x_n$ with probabilities
  $\frac{1}{2}$ each.
\end{theorem}
\begin{proof}
  See \citet{discvar}.
\end{proof}
In both cases we obtain that the maximum variance is achieved for $p_{i(1)} =
p_{i(-1)} = \frac{1}{2}$ and is equal to $1$, so $\sigma_{ii} \in 
[0, 1]$  and $|\sigma_{ij}| \in [0, 1]$. Furthermore, we can also prove that
the eigenvalues of $\Sigma$ are bounded using the same arguments as in Lemma
\ref{thm:mvebereigen}.
\begin{lemma}
\label{thm:dirlambda}
  Let $\mathbf{T} \sim Tri_k(\mathbf{p})$, and let $\Sigma$ be its covariance
  matrix. Let $\lambda_i$, $i = 1, \ldots, k$ be the eigenvalues of $\Sigma$.
  Then
  \begin{align*}
    &0 \leqslant \sum_{i=1}^k \lambda_i \leqslant k&
    &\text{and}&
    &0 \leqslant \lambda_i \leqslant k.
  \end{align*}
\end{lemma}
\begin{proof}
  See the proof of Lemma \ref{thm:mvebereigen} in Appendix \ref{app:proofs}.
\end{proof}
These bounds define again a closed convex set in $\mathbb{R}^k$, described by
the family
\begin{equation*}
  \mathcal{L} = \left\{ \Delta^{k-1}(c) : c \in \left[ 0, k \right]\right\},
\end{equation*}
where $\Delta^{k-1}(c)$ is the non-standard $k-1$ simplex from Equation \ref{eq:simplex}.

Another useful result, which we will use in Section \ref{sec:dagprop} to link
inference on UGs and DAGs, is introduced below.
\begin{theorem}
\label{thm:triber2}
  Let $\mathbf{T} \sim Tri_k(\mathbf{p})$; then $|\mathbf{T}| = \mathbf{B}
  \sim Ber_k(\mathbf{p}^*)$ and \\ $|T_i| = B_i \sim Ber(p^*)$.
\end{theorem}
\begin{proof}
See Appendix \ref{app:proofs}.
\end{proof}
It follows that the variance of each $T_i$ can be decomposed in two parts:
\begin{equation}
\label{eqn:vardecomp}
  \VAR(T_i) = \VAR(B_i) + 4 p_{i(1)} p_{i(-1)}.
\end{equation}
The first is a function of the corresponding component $|T_i| = B_i$ of
the transformed random vector, while the second depends only on the probabilities
associated with $-1$ and $1$ (which correspond to $\overleftarrow{a_{ij}}$ and
$\overrightarrow{a_{ij}}$ in Equation \ref{eqn:tridef}).

\pagebreak

\section{Properties of $\eprior$ and $\eposterior$}
\label{sec:properties}

The results derived in the previous section provide the foundation for characterising
$\eprior$ and $\eposterior$. To this end, it is useful to distinguish three
cases corresponding to different configurations of the probability mass
among the graph structures $\mathcal{G}(\mathcal{E}) \in \mathbf{G}$:
\begin{itemize}
  \item \textit{minimum entropy}: the probability mass is concentrated on a
    single graph structure. This is the best possible configuration for
    $\eposterior$, because only one edge set $E$ (or one arc set $A$) has
    a non-zero posterior probability. In other words, the data $\mathcal{D}$
    provide enough information to identify a single graph $\mathcal{G}$ with
    posterior probability $1$;
  \item \textit{intermediate entropy}: several graph structures have non-zero
    probabilities. This is the case for informative priors $\eprior$ and for
    the posteriors $\eposterior$ resulting from real-world data sets;
  \item \textit{maximum entropy}: all graph structures in $\mathbf{G}$ have the
    same probability. This is the worst possible configuration for $\eposterior$,
    because it corresponds to the non-informative prior from Equation
    \ref{eqn:flatprior}. In other words, the data $\mathcal{D}$ do not provide
    any information useful in identifying a high-posterior graph $\mathcal{G}$.
\end{itemize}
Clearly, \textit{minimum} and \textit{maximum entropy} are limiting cases for
$\eposterior$; the former is non-informative about $\mathcal{G} (\mathcal{E})$,
while the latter identifies a single graph in $\mathbf{G}$. As we will show in
Sections \ref{sec:ugprop} (for UGs) and \ref{sec:dagprop} (for DAGs), they
provide useful reference points in determining which edges (or arcs) have
significant posterior probabilities and in analysing the variability of the
graph structure.

\subsection{Undirected graphs}
\label{sec:ugprop}

In the \textit{minimum entropy} case, only one configuration of edges $E$ has
non-zero probability, which means that
\begin{align*}
  &p_{ij} = \left\{
    \begin{aligned}
      &1& &\text{if $e_{ij} \in E$}     \\
      &0& &\text{otherwise}&
    \end{aligned}
    \right.&
  &\text{and}&
  &\Sigma = \mathbf{O}.
\end{align*}

The uniform distribution over $\mathbf{G}$ arising from the \textit{maximum
entropy} case has been studied extensively in random graph theory \citep{bollobas};
its two most relevant properties are that all edges $e_{ij}$ are independent
and have $p_{ij} = \frac{1}{2}$. As a result, $\Sigma = \frac{1}{4}I_k$; all
edges display their maximum possible variability, which along with the fact
that they are independent makes this distribution non-informative for
$\mathcal{E}$ as well as $\mathcal{G}(\mathcal{E})$.

The \textit{intermediate entropy} case displays a middle-ground behaviour between
the \textit{minimum} and \textit{maximum entropy} cases. The expected value and
the covariance matrix of $\mathcal{E}$ do not have a definite form beyond the
bounds derived in Section \ref{sec:mvber}. When considering posteriors arising
from  real-world data, we have in practice that most edges in $\mathcal{E}$
represent conditional dependence relationships that are completely unsupported
by the data. This behaviour has been explained by \citet{causality} with the
tendency of ``good'' graphical models to represent the causal relationships
underlying the data, which are typically sparse.
As a result, we have that $\E(e_{ij}) = 0$ and $\VAR(e_{ij}) = 0$ for many
$e_{ij}$, so $\Sigma$ is almost surely singular unless such edges are excluded
from the analysis. Edges that appear with $p_{ij} \simeq \frac{1}{2}$ have about
the same marginal probability and variance as in the \textit{maximum entropy}
case, so their marginal behaviour is very close to random noise. On the other
hand, edges with probabilities near $0$ or $1$ can be considered to have a good
support (against or in favour, respectively). As $p_{ij}$ approaches $0$ or $1$,
$e_{ij}$ approaches its \textit{minimum entropy}.

The closeness of a multivariate Bernoulli distribution to the \textit{minimum}
and \textit{maximum entropy} cases can be represented in an intuitive way by
considering the eigenvalues $\boldsymbol{\lambda} = [\lambda_1, \ldots,
\lambda_k]^T$ of its covariance matrix $\Sigma$. Recall that the
$\boldsymbol{\lambda}$ can assume values in the convex set $\mathcal{L}$ defined
in Equation \ref{eq:simplex}, which corresponds to the region of the first
orthant delimited by the non-standard simplex $\Delta^{k-1}(\frac{k}{4})$. In
the \textit{minimum entropy} case we have that $\Sigma = \mathbf{O}$, so
$\lambda_1 = \ldots = \lambda_k = 0$, and in the \textit{maximum entropy case}
$\Sigma = \frac{1}{4}I_k$, so $\lambda_1 = \ldots = \lambda_k = \frac{1}{4}$;
both points lie on the boundary of $\mathcal{L}$, the first in the origin and
the second in the middle of $\Delta^{k-1}(\frac{k}{4})$. The distance between
$\boldsymbol{\lambda}$ and these two points provides an intuitive way of measuring
the variability of $\mathcal{E}$ and, indirectly, the entropy of the corresponding
probability distributions $\eposterior$ and $\eprior$. It is important to note,
however, that different distributions over $\mathbf{G}$ may have identical first
and second order moments when modelled through $\mathcal{E}$. Such distributions
will have the same $\boldsymbol{\lambda}$ and will therefore map to the same
point in $\mathcal{L}$.

A simple example comprising three different distributions over a set of two edges
is illustrated below.

\begin{example}
\label{ex:base}
  Consider three multivariate Bernoulli distributions $\mathbf{B}_1$, $\mathbf{B}_2$,
  $\mathbf{B}_3$ over two edges (denoted with $e_1 \sim E_1$ and $e_2 \sim E_2$ for
  brevity) with covariance matrices
  \begin{align*}
    &\Sigma_1 = \begin{bmatrix} 0.24 & 0.04 \\ 0.04 & 0.24 \end{bmatrix},&
    &\Sigma_2 = \begin{bmatrix} 0.1056 & -0.0336 \\ -0.0336 & 0.2016 \end{bmatrix},&
    &\Sigma_3 = \begin{bmatrix} 0.1056 & 0.1456 \\ 0.1456 & 0.2016 \end{bmatrix}
  \end{align*}
  and eigenvalues
  \begin{align*}
    &\boldsymbol{\lambda}_1 = \begin{bmatrix} 0.28 \\ 0.20 \end{bmatrix},&
    &\boldsymbol{\lambda}_2 = \begin{bmatrix} 0.2121 \\ 0.095 \end{bmatrix},&
    &\boldsymbol{\lambda}_3 = \begin{bmatrix} 0.3069 \\ 0.0003 \end{bmatrix}.
  \end{align*}
  Their positions in $\mathcal{L}$ are shown in Figure \ref{fig:base}. $\mathbf{B}_1$
  is the closest to $\left(\frac{1}{4}, \frac{1}{4}\right)$, the point corresponding
  to the maximum entropy case, while $\mathbf{B}_2$ and $\mathbf{B}_3$ are farther
  from $\left(\frac{1}{4}, \frac{1}{4}\right)$ than $\mathbf{B}_1$ due to the
  increasing correlation between $e_1$ and $e_2$ (which are independent in the
  \textit{maximum entropy} case). The correlation coefficients for $\mathbf{B}_1$,
  $\mathbf{B}_2$ and $\mathbf{B}_3$ are $\COR_{\mathbf{B}_1}(E_1, E_2) = 0.1666$,
  $\COR_{\mathbf{B}_2}(E_1, E_2) = -0.2303$, $\COR_{\mathbf{B}_3}(E_1, E_2) =
  0.9978$, and they account for the increasing difference between the eigenvalues
  of each covariance matrix. In fact, $\Sigma_3$ is nearly singular because of the
  strong linear relationship between $e_1$ and $e_2$, and it is therefore very
  close to one of the axes delimiting the first quadrant.

  \begin{figure}[t]
    \begin{center}
      \includegraphics{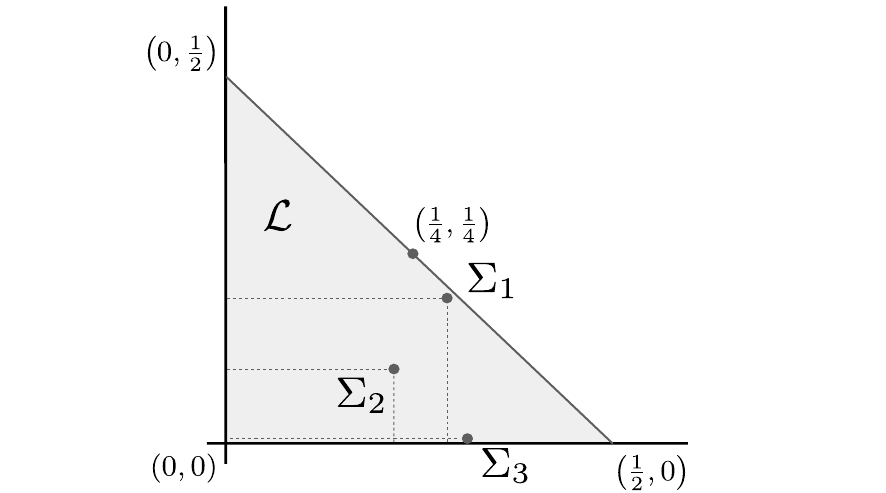}
    \end{center}
    \caption{The covariance matrices $\Sigma_1$, $\Sigma_2$ and $\Sigma_3$ from
      Example \ref{ex:base} represented as functions of their eigenvalues in the
      convex set $\mathcal{L}$. The points $(0,0)$ and $(\frac{1}{4},
      \frac{1}{4})$ correspond to the \textit{minimum entropy} and \textit{maximum
      entropy} cases.}
  \label{fig:base}
  \end{figure}

  If we denote with $E_{00} = \{\varnothing\}$, $E_{01} = \{e_2\}$, $E_{10} = \{e_1\}$,
  and $E_{11} = \{e_1, e_2\}$ all possible edge sets and with $p_{00}$, $p_{01}$,
  $p_{10}$ and $p_{11}$ the associated probabilities, for $\mathbf{B}_1$ we have
  \begin{align*}
    &p_{00} = 0.2,& &p_{01} = 0.2,& &p_{10} = 0.2& &\text{and}& &p_{11} = 0.4.
  \end{align*}
  This is indeed close to a uniform distribution. The probability of both
  $e_1$ and $e_2$ is $0.6$ and the variance is $0.24$, which are again similar
  to the reference values for the \textit{maximum entropy} case. On the other
  hand, for $\mathbf{B}_2$ we have
  \begin{align*}
    &p_{00} = 0,& &p_{01} = 0.12,& &p_{10} = 0.28& &\text{and}& &p_{11} = 0.6.
  \end{align*}
  These probabilities are markedly different from a uniform distribution; the
  probabilities of $e_1$ and $e_2$ are respectively $0.88$ and $0.72$. Considering
  also the correlation between $e_1$ and $e_2$, it is intuitively clear why
  $\Sigma_2$ is not as close as $\Sigma_1$ to $\left(\frac{1}{4}, \frac{1}{4}\right)$.
  This is also true for $\mathbf{B}_3$, which has the same marginal distributions
  as $\mathbf{B}_2$ but with a much stronger correlation.
\end{example}

\subsection{Directed acyclic graphs}
\label{sec:dagprop}

The behaviour of the multivariate Trinomial distribution in the \textit{minimum}
and \textit{intermediate entropy} cases is similar to the one of the multivariate
Bernoulli in many respects, but presents profound differences in the \textit{maximum
entropy} case. The reason for these differences is that the structure of a Bayesian
network is assumed to be acyclic. Therefore, the state of each arc (i.e. whether is
present in the DAG and its direction) is influenced by the state of all other
possible arcs even in the \textit{maximum entropy} case, when otherwise they would
be independent. Furthermore, the acyclicity constraint cannot be written in closed
form, making the derivation of exact results on the moments of the distribution of
$\mathcal{E}$ particularly difficult.

To obtain some simple expressions for the expected value and the covariance matrix,
we will first prove a simple theorem on DAGs, which essentially states that if we
reverse the direction of every arc the resulting graph is still a DAG.
\begin{theorem}
\label{thm:acyclic}
  Let $G = (\mathbf{V}, A)$ be a DAG, and let $G^* = (\mathbf{V}, A^*)$ another
  directed graph such that
  \begin{align*}
    &\overrightarrow{a_{ij}} \in A^* \Longleftrightarrow \overleftarrow{a_{ij}} \in A&
    &\text{and}&
    &\overleftarrow{a_{ij}} \in A^* \Longleftrightarrow \overrightarrow{a_{ij}} \in A
  \end{align*}
  for every $a_{ij} \in A$. Then $G^*$ is also acyclic.
\end{theorem}
\begin{proof}
  See Appendix \ref{app:proofs}.
\end{proof}
An immediate consequence of this theorem is that for every DAG including the
arc $\overrightarrow{a_{ij}}$ there exists another DAG including the arc
$\overleftarrow{a_{ij}}$. Since all DAGs have the same probability in the
\textit{maximum entropy case}, this implies that both directions of every arc
have the same probability,
\begin{align}
\label{eqn:1storder}
  &\overrightarrow{p_{ij}} = \overleftarrow{p_{ij}}& &\text{for every possible } a_{ij}, i \neq j.
\end{align}
Then the expected value of each marginal Trinomial distribution is equal to
\begin{equation*}
  \E(A_{ij}) = \overrightarrow{p_{ij}} - \overleftarrow{p_{ij}} = 0
\end{equation*}
and its variance is equal to
\begin{equation*}
  \VAR(A_{ij}) = \overrightarrow{p_{ij}} + \overleftarrow{p_{ij}} - (\overrightarrow{p_{ij}} - \overleftarrow{p_{ij}})^2
    = 2\overrightarrow{p_{ij}}.
\end{equation*}

The joint probabilities associated with each pair of arcs also symmetric in the
maximum entropy case, again due to Theorem \ref{thm:acyclic}. Denote with 
$\mathring{a_{ij}}$ the event that arc $a_{ij}$ is not present in the DAG. If
we consider that both directions of every arc have the same probability and that
there is no explicit ordering among the arcs, we have
\begin{gather}
\label{eqn:2ndorder}
\begin{split}
  \Prob(\overrightarrow{a_{ij}}, \overrightarrow{a_{kl}}) = \Prob(\overleftarrow{a_{ij}}, \overleftarrow{a_{kl}}),
  \qquad \qquad
  \Prob(\overrightarrow{a_{ij}}, \overleftarrow{a_{kl}}) = \Prob(\overleftarrow{a_{ij}}, \overrightarrow{a_{kl}}), \\
  \Prob(\mathring{a_{ij}}, \overrightarrow{a_{kl}}) = \Prob(\overrightarrow{a_{ij}}, \mathring{a_{kl}}) =
  \Prob(\mathring{a_{ij}}, \overleftarrow{a_{kl}}) = \Prob(\overleftarrow{a_{ij}}, \mathring{a_{kl}}). 
\end{split}
\end{gather}
Then the expression for the covariance simplifies to
\begin{equation*}
  \COV(A_{ij}, A_{kl}) = 2 \left[ \Prob(\overrightarrow{a_{ij}}, \overrightarrow{a_{kl}}) - \Prob(\overrightarrow{a_{ij}}, \overleftarrow{a_{kl}}) \right],
\end{equation*}
which can be interpreted as the difference in probability between a \textit{serial
connection} (i.e. $v_i \rightarrow v_j \rightarrow v_l$, if $j = k$) and a
\textit{converging connection} (i.e. $v_i \rightarrow v_j \leftarrow v_l$) if
the arcs are incident on a common node \citep{jensen}. This is interesting
because v-structures are invariant within equivalence classes, while other
patterns of arcs are not \citep{chickering}; indeed, equivalence classes are
usually represented as \textit{partially directed acyclic graphs} (PDAGs) in
which only arcs belonging to v-structures are directed. All other arcs, with
the exclusion of those which could introduce additional v-structures (known
as \textit{compelled arcs}), are replaced with the corresponding (undirected)
edges. Therefore, the combination of high values of $|\COV(A_{ij}, A_{kl})|$ 
and $\mathring{p_{ij}}$ is indicative of the belief that the corresponding arcs
are directed in the PDAG identified by the equivalence class. Along with with
$\VAR(A_{ij})$ and $\VAR(A_{kl})$, it is also indicative of the stability of
the graph structure, both in the arcs and their directions. In an uninformative
prior, such as the distribution we are now considering in the \textit{maximum 
entropy} case, we expect all covariances to be small; we will show this is the
case in Theorem \ref{thm:hoeffding}. On the other hand, in an informative 
distribution such as the ones considered in the \textit{intermediate entropy}
case, we expect covariances to be closer to their upper bounds for arcs that
are compelled or part of a converging connection, and closer to zero for arcs
whose direction is not determined in the equivalence class. Note that the sign
of $\COV(A_{ij}, A_{kl})$ depends on the way the two possible directions of
each arc are associated with $1$ and $-1$; a simple way to obtain a consistent
parameterisation is to follow the natural ordering of the variables (i.e. if
$i \leqslant j$ then the arc incident on these nodes is taken to be $A_{ij}$,
$\overrightarrow{a_{ij}}$ is associated with $1$ and $\overleftarrow{a_{ij}}$
with $-1$).

\begin{figure}[t]
  \begin{center}
  \includegraphics[width=0.9\textwidth]{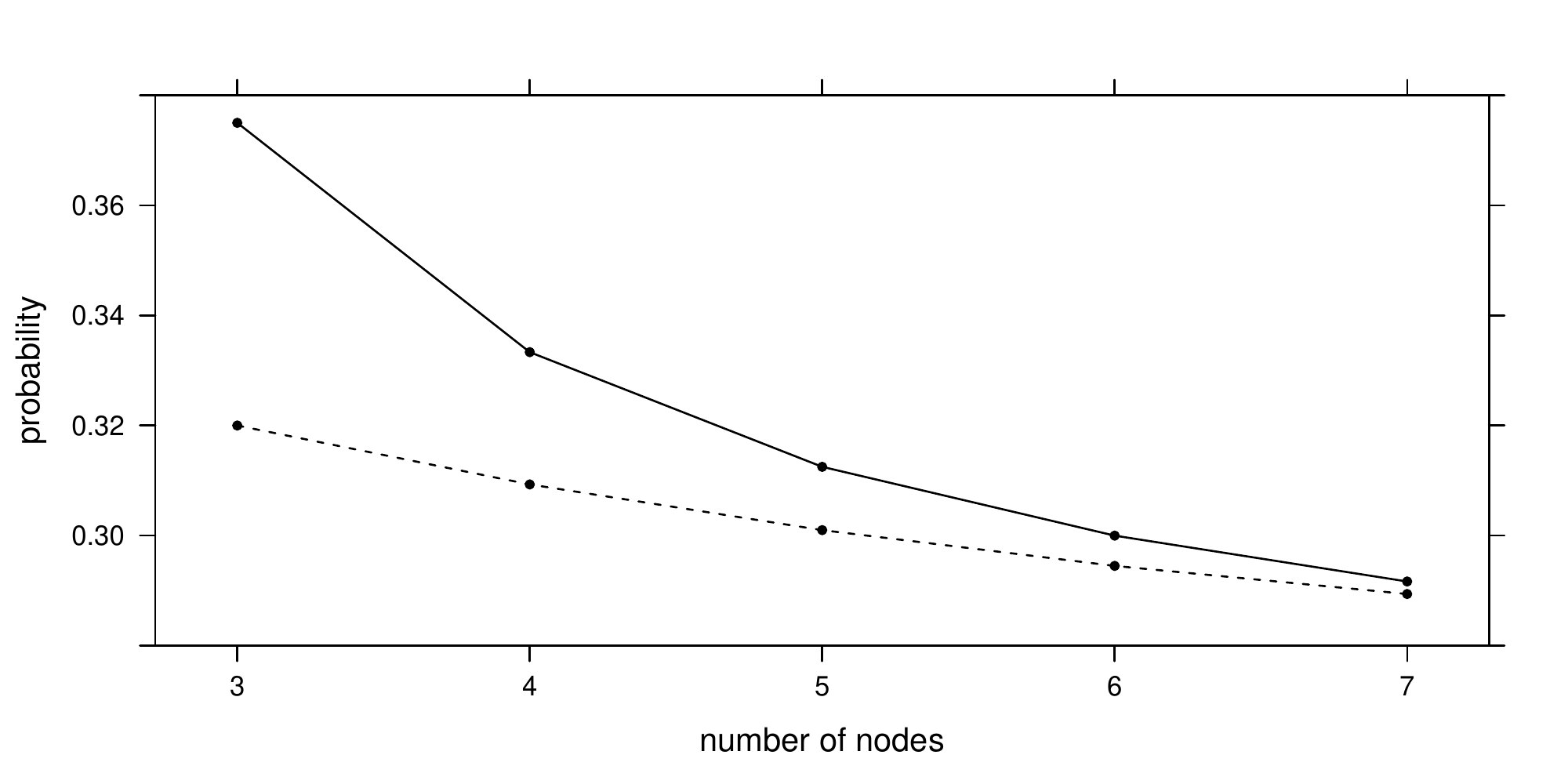}
  \caption{Exact (dashed line) and approximate (solid line) probabilities of an
    arc being present in a DAG with $3$, $4$, $5$, $6$, and $7$ nodes. The dotted
    line represents the limiting value in the number of nodes.}
  \label{fig:exactarcprob}
  \end{center}
\end{figure}
\begin{figure}[t]
  \begin{center}
  \includegraphics[width=0.9\textwidth]{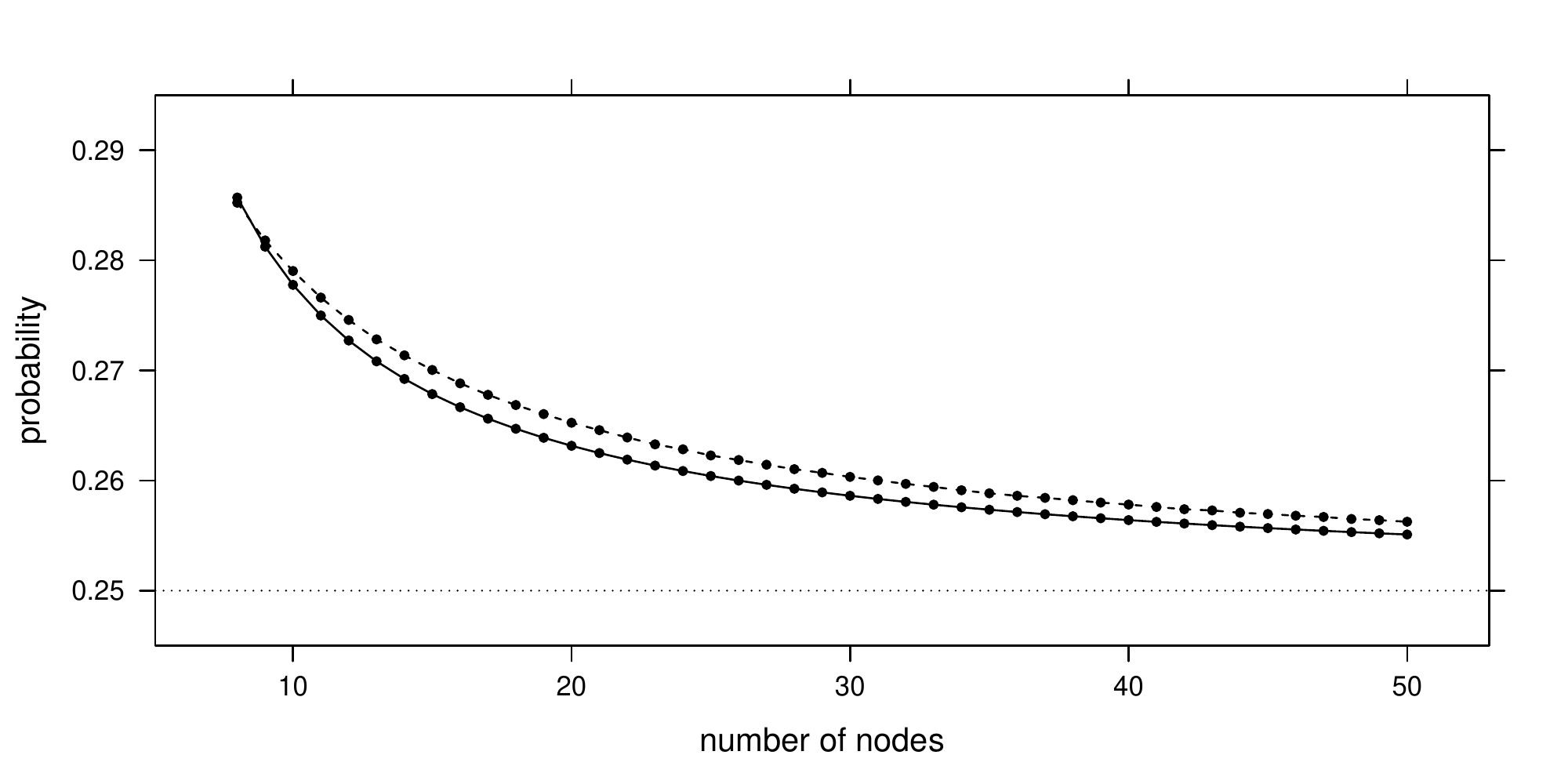}
  \caption{Estimated (dashed line) and approximate (solid line) probabilities of
    an arc being present in a DAG with $8$ to $50$ nodes. The dotted line
    represents the limiting value in the number of nodes.}
  \label{fig:mtarcprob}
  \end{center}
\end{figure}

The equalities in Equations \ref{eqn:1storder} and \ref{eqn:2ndorder} drastically
reduce the number of free parameters in the \textit{maximum entropy} case. The
marginal distribution of each arc now depends only on $\overrightarrow{p_{ij}}$,
whose value can be derived from the following numerical approximation by
\citet{melancon2}.

\begin{theorem}
\label{thm:idecozman}
  The average number of arcs in a DAG with $n$ nodes is approximately $\frac{1}{4}n^2$
  in the \textit{maximum entropy} case.
\end{theorem}
\begin{proof}
  See \citet{melancon2}.
\end{proof}
\begin{theorem}
\label{thm:corollary}
  Let $G = (\mathbf{V}, A)$ be a DAG with $n$ nodes. Then for each possible arc
  $a_{ij}, i \neq j$ we have that in the maximum entropy case
  \begin{align*}
    &\overrightarrow{p_{ij}} = \overleftarrow{p_{ij}} \simeq \frac{1}{4} + \frac{1}{4(n-1)}&
    &\text{and}&
    &\mathring{p_{ij}} \simeq \frac{1}{2} - \frac{1}{2(n-1)}.
  \end{align*}
\end{theorem}
\begin{proof}
  See Appendix \ref{app:proofs}.
\end{proof}

The quality of this approximation is examined in Figure \ref{fig:exactarcprob}
and Figure \ref{fig:mtarcprob}. In Figure \ref{fig:exactarcprob}, the values
provided by Theorem \ref{thm:corollary} for DAGs with $3$, $4$, $5$, $6$ and
$7$ nodes are compared to the corresponding true values. The latter have been
computed by enumerating all possible DAGs of that size (i.e. the whole
population) and computing the relative frequency of each possible arc. In Figure
\ref{fig:mtarcprob}, the values provided by Theorem \ref{thm:corollary} for
DAGs with $8$ to $50$ nodes are compared with the corresponding estimated values
computed over a set of $10^9$ DAGs of the same size. The latter have been
generated with uniform probability using the algorithm from \citet{melancon}
as implemented in the bnlearn package \citep{jss09,bnlearn} for R \citep{R}.

We can clearly see that the approximate values are close to the corresponding
true (in Figure \ref{fig:exactarcprob}) or estimated (in Figure \ref{fig:mtarcprob})
values for DAGs with at least $6$ nodes. This is not a significant limitation;
the true values can be easily computed via exhaustive enumeration for DAGs with
$3$, $4$ and $5$ nodes (they are reported in Appendix \ref{app:numbers}, along
with other relevant quantities). Furthermore, it is evident both from Theorem
\ref{thm:corollary} and from Figures \ref{fig:exactarcprob} and \ref{fig:mtarcprob}
that, as the number of nodes diverges,
\begin{align}
\label{eqn:limits}
  &\lim_{n \to \infty} \overrightarrow{p_{ij}} = \lim_{n \to \infty} \overleftarrow{p_{ij}} = \frac{1}{4}&
  &\text{and}&
  &\lim_{n \to \infty} \mathring{p_{ij}} = \frac{1}{2}.
\end{align}
If we take the absolute value of this asymptotic Trinomial distribution, the
resulting random variable is $Ber(p_{ij})$ with $p_{ij} = \frac{1}{2}$, which
is the marginal distribution of an edge in an UG in the \textit{maximum entropy}
case. The absolute value transformation can be interpreted as ignoring the
direction of the arc; the events $\overleftarrow{a_{ij}} \in A$ and
$\overrightarrow{a_{ij}} \in A$ collapse into $e_{ij} \in E$, while
$\overleftarrow{a_{ij}}, \overrightarrow{a_{ij}} \notin A$ maps to $e_{ij} 
\notin E$. As a result, the marginal distribution of an arc is remarkably
similar to the one of the corresponding edge in an undirected graph for
sufficiently large DAGs; in both cases, the nodes $v_i$ and $v_j$ are linked
with probability $\frac{1}{2}$.

No result similar to Theorem \ref{thm:idecozman} has been proved for arbitrary
pairs of arcs in a directed acyclic graph; therefore, the structure of the
covariance matrix $\Sigma$ can be derived only in part. Variances can be
approximated using the approximate probabilities from Theorem \ref{thm:corollary}:
\begin{equation}
  \label{eqn:approxvar}
  \VAR(A_{ij}) = 2\overrightarrow{p_{ij}} \simeq \frac{1}{2} + \frac{1}{2(n-1)} \to \frac{1}{2} \,\text{ as }\, n \to \infty.
\end{equation}

\begin{figure}[p]
  \begin{center}
  \includegraphics[width=0.9\textwidth]{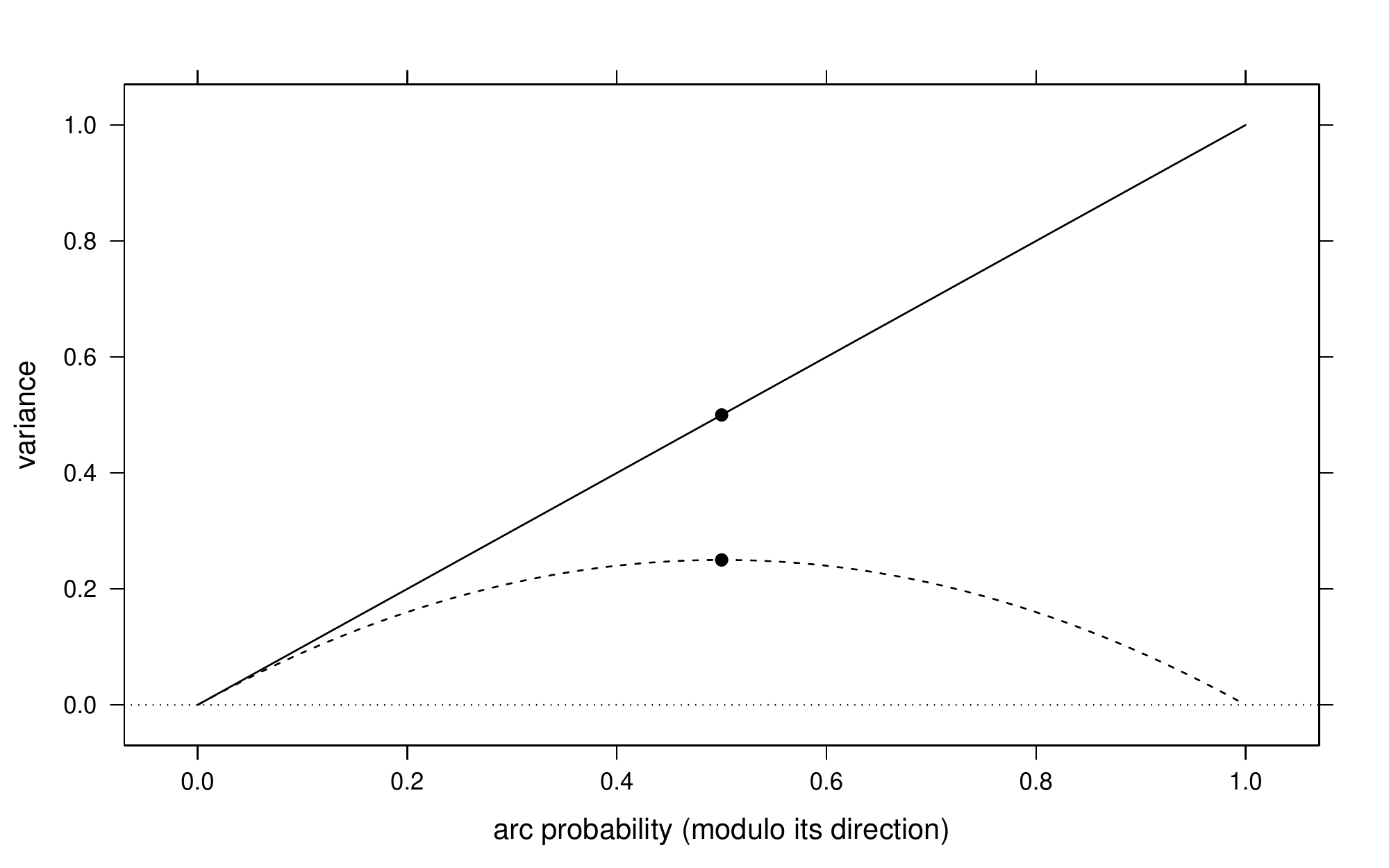}
  \caption{Decomposition of the asymptotic variance of an arc in the part that
    depends only on its presence (dashed line) and the part that depends only
    on its direction (solid line). The dots correspond to the respective values
    in the maximum entropy case.}
  \label{fig:vardecomp}
  \includegraphics[width=0.9\textwidth]{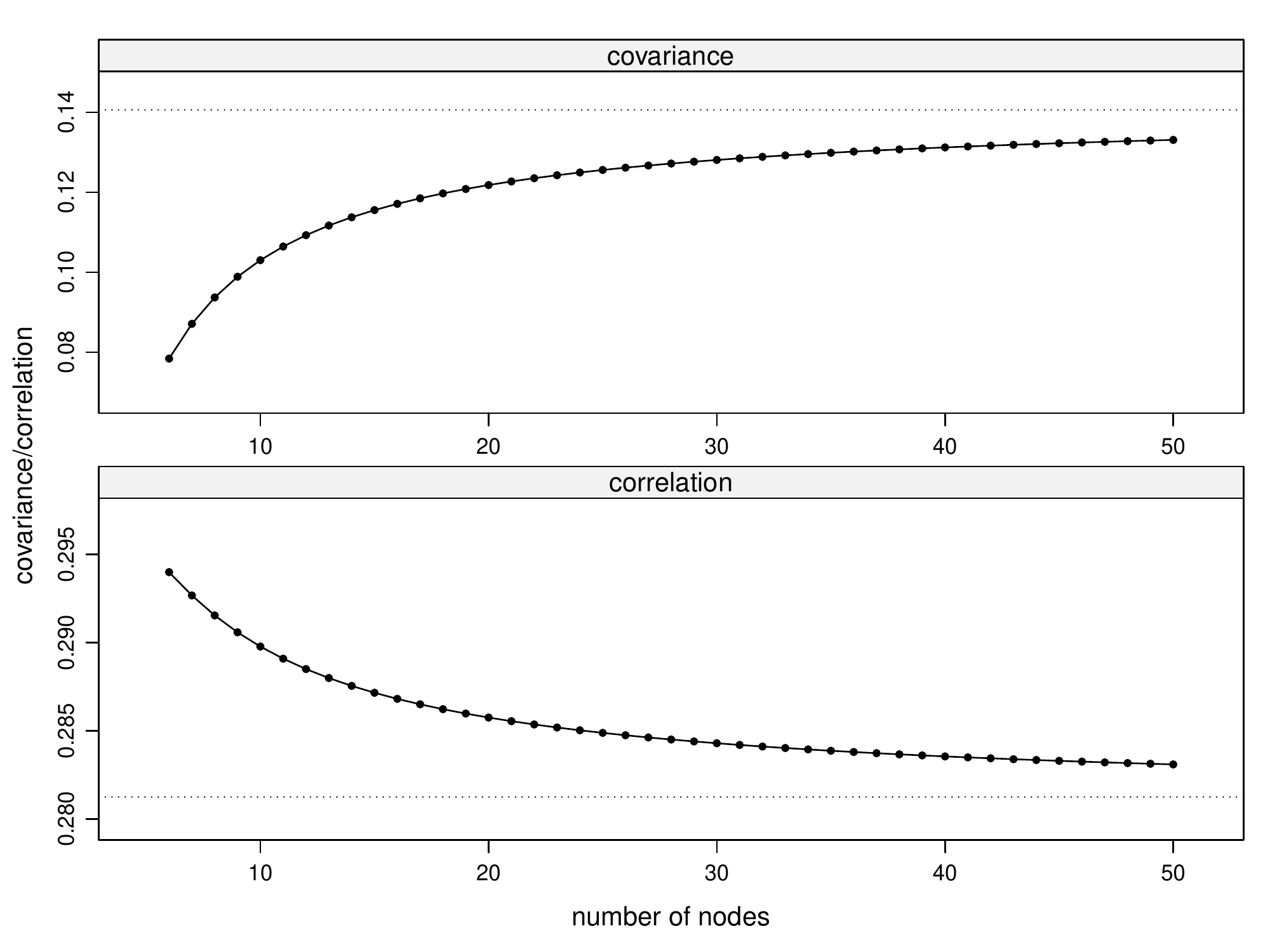}
  \caption{Bounds for the absolute value of the covariance and the correlation
     coefficient of two arcs in a DAG with $6$ to $50$ nodes. The dotted lines
     represent the respective limiting values.}
  \label{fig:approxcor}
  \end{center}
\end{figure}

Therefore, maximum variance (of each arc) and maximum entropy (of the graph
structure) are distinct, as opposed to what happens in UGs. However, we can use
the decomposition of the variance introduced in Equation \ref{eqn:vardecomp}
to motivate why the \textit{maximum entropy} case is still a ``worst case''
outcome for $\eposterior$. As we can see from Figure \ref{fig:vardecomp}, the
contributions of the presence of an arc (given by the transformation $|A_{ij}|$)
and its direction (given by the $4\overrightarrow{p_{ij}}\overleftarrow{p_{ij}}
= 4\overrightarrow{p_{ij}}^2$ term) to the variance are asymptotically equal.
This is a consequence of the limits in Equation \ref{eqn:limits}, which imply
that an arc (modulo its direction) has the same probability to be present in
or absent from the DAG and that its directions also have the same probability.
As a result, we are not able to make any decision about either the presence
of the arc or its direction. On the contrary, when $\VAR(A_{ij})$ reaches it
maximum at $1$ we have that $\Prob(\{\overrightarrow{a_{ij}}, 
\overleftarrow{a_{ij}}\}) = 1$ and $\Prob(\mathring{a_{ij}}) = 0$, so we are
sure that the arc will be present in the DAG in one of its two possible directions.

As for the covariances, it is possible to obtain tight bounds using
\textit{Hoeffding's identity} \citep{hoeffding,hoeffding2},
\begin{equation}
\label{eqn:hoeffding}
  \COV(X, Y) = \iint_{\mathbb{R}^2} F_{X, Y}(x, y) - F_X(x)F_Y(y) dxdy,
\end{equation}
and the decomposition of the joint distribution of dependent random variables
provided by the \textit{Farlie-Morgenstern-Gumbel} (FMG) family of distributions
\citep{fmg}, which has the form
\begin{align}
\label{eqn:fmg}
  &F_{X, Y}(x, y) = F_X(x)F_Y(y) \left[ 1 + \varepsilon (1 - F_X(x))(1 - F_Y(y))\right],&
  &|\varepsilon| \leqslant 1.
\end{align}
In Equations \ref{eqn:hoeffding} and \ref{eqn:fmg}, $F_{X,Y}$, $F_X$ and $F_Y$
denote the cumulative distribution functions of the joint and marginal
distributions of $X$ and $Y$, respectively.

\begin{theorem}
\label{thm:hoeffding}
  Let $G = (V, A)$ be a DAG, and let $a_{ij}$, $i \neq j$ and $a_{kl}$, $k \neq l$
  be two possible arcs. Then in the \textit{maximum entropy} case we have that
  \begin{equation}
    \label{eqn:boundcov}
    \left|\COV(A_{ij}, A_{kl})\right| \lessapprox 4\left[\frac{3}{4} - \frac{1}{4(n -1)}\right]^2 \left[\frac{1}{4} + \frac{1}{4(n - 1)}\right]^2
  \end{equation}
  and
  \begin{equation}
    \label{eqn:boundcor}
    \left|\COR(A_{ij}, A_{kl})\right| \lessapprox 2\left[\frac{3}{4} - \frac{1}{4(n -1)}\right]^2 \left[\frac{1}{4} + \frac{1}{4(n - 1)}\right].
  \end{equation}
\end{theorem}
\begin{proof}
  See Appendix \ref{app:proofs}.
\end{proof}

The bounds obtained from this theorem appear to be tight in the light of the
true values for the covariance and correlation coefficients (computed again by
enumerating all possible DAGs of size $3$ to $7$). Figure \ref{fig:approxcor}
shows the bounds for DAGs with $6$ to $50$ nodes; for DAGs with $3$, $4$ and
$5$ nodes the approximation of $\overrightarrow{p_{ij}}$ the bounds are based
on is loose, and the true values of covariance and correlation are known. Non-null
covariances range from $\pm 0.08$ (for DAGs with $3$ nodes) to $\pm 0.08410$
(for DAGs with $7$ nodes), while non-null correlation coefficients vary from
$\pm 0.125$ (for DAGs with $3$ nodes) to $\pm 0.1423$ (for DAGs with $7$ nodes).
Both covariance and correlation appear to be strictly increasing in modulus as
the number of nodes increases, and converge to the limiting values of the bounds
($0.140625$ and $0.28125$, respectively) from below.

\begin{figure}[t]
  \begin{center}
  \includegraphics[width=0.9\textwidth]{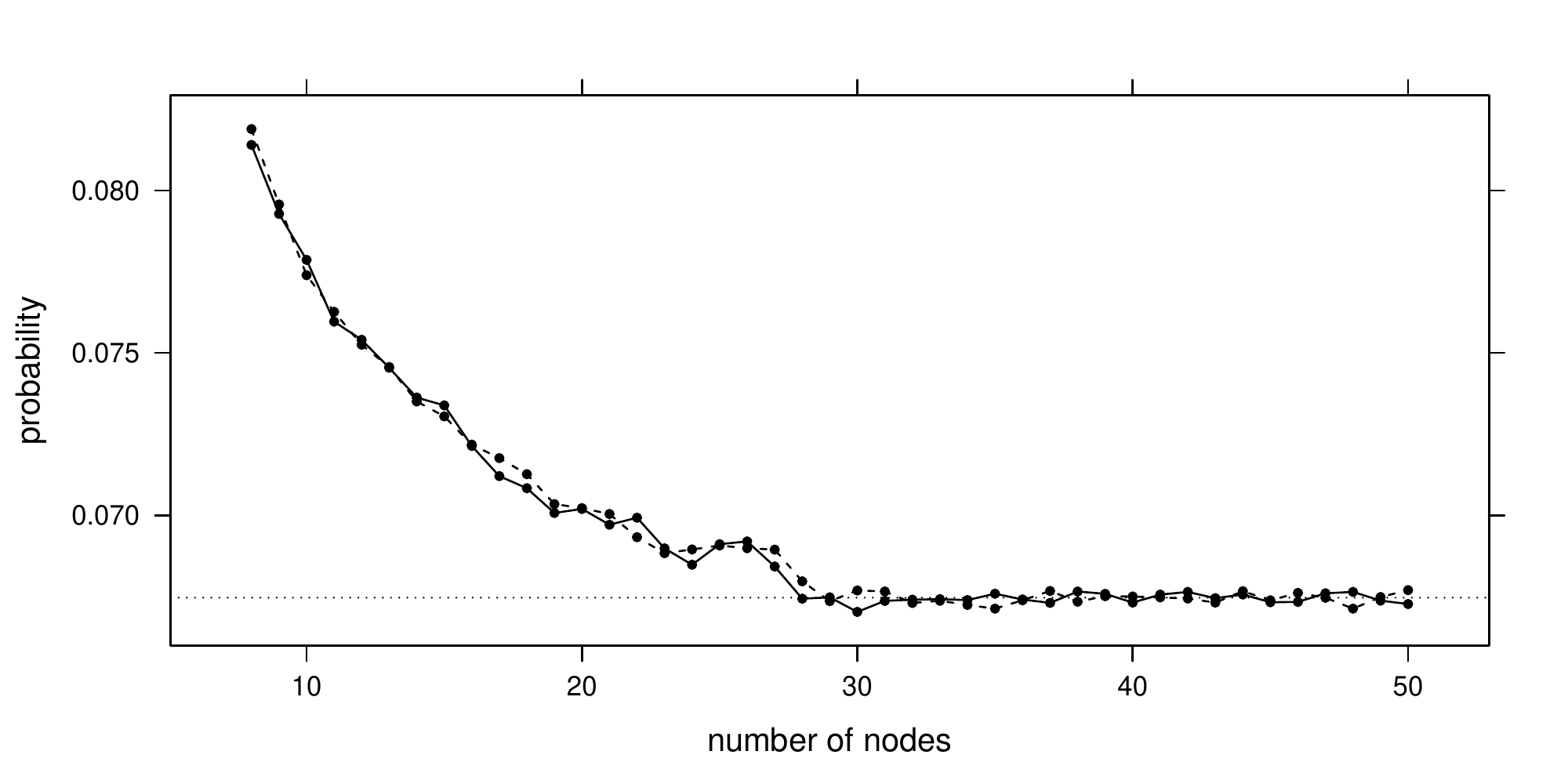}
  \caption{Approximate values for $\Prob(\protect\overrightarrow{a_{ij}},
     \protect\overrightarrow{a_{kl}})$ (solid line) and
     $\Prob(\protect\overrightarrow{a_{ij}}, \protect\overleftarrow{a_{kl}})$
     (dashed line) for DAGs with $8$ to $50$ nodes. The dotted line represents
     their asymptotic value.}
  \label{fig:approxp}
  \end{center}
\end{figure}

Some other interesting properties are apparent from true values of the covariance
coefficients reported in Appendix \ref{app:numbers}. They are reported below as
conjectures because, while they describe a systematic behaviour that emerges from
the DAGs whose sizes we have a complete enumeration for, we were not able to
substantiate them with formal proofs.
\begin{conjecture}
\label{conj:unc}
  Arcs that are not incident on a common node are uncorrelated.
\end{conjecture}
This is a consequence of the fact that if we consider $A_{ij}$ and $A_{kl}$ with
$i \neq j \neq k \neq l$, we have $\Prob(\overrightarrow{a_{ij}}, \overrightarrow{a_{kl}})
= \Prob(\overrightarrow{a_{ij}}, \overleftarrow{a_{kl}})$. Therefore
$\COV(A_{ij}, A_{kl}) = 0$. This property seems to generalise to DAGs with
more than $7$ nodes. Figure \ref{fig:approxp} shows approximate estimates for 
$\Prob(\overrightarrow{a_{ij}}, \overrightarrow{a_{kl}})$ and 
$\Prob(\overrightarrow{a_{ij}}, \overleftarrow{a_{kl}})$ for DAGs with $8$ to
$50$ nodes, obtained again from $10^9$ DAGs generated with uniform probability.
The curves for the two probabilities are overlapping and very close to each
other for all the considered DAG sizes, thus supporting Conjecture \ref{conj:unc}.

\begin{conjecture}
  The covariance matrix $\Sigma$ is sparse.
\end{conjecture}
The proportion of arcs incident on a common node converges to zero as the number of
nodes increases; therefore, if we assume Conjecture \ref{conj:unc} is true, the
proportion of elements of $\Sigma$ that are equal to $0$ has limit
\begin{equation}
  1 \geqslant \lim_{n \to \infty} \frac{{n \choose 2} {n-2 \choose 2}}{{n \choose 2}{n \choose 2} - {n \choose 2}}
    \geqslant \lim_{n \to \infty} \frac{(n -2)(n - 3)}{n(n-1)} = 1.
\end{equation}
Furthermore, even arcs that are incident on a common node are not strongly correlated.
\begin{conjecture}
\label{conj:increasing}
  Both covariance and correlation between two arcs incident on a common node
  are monotonically increasing in modulus.
\end{conjecture}
\begin{conjecture}
  The covariance between two arcs incident on a common node takes values in the
  interval $[0.08,$ $0.140625]$ in modulus, while the correlation takes values
  in $[0.125, 0.28125]$ in modulus.
\end{conjecture}
These intervals can be further reduced to $[0.08410, 0.140625]$ and $[0.1423,$
$0.28125]$ for DAGs larger than $7$ nodes due to Conjecture \ref{conj:increasing}.

As far as the other two cases are concerned, in the \textit{minimum entropy} case
we have that
\begin{align*}
  &\E(A_{ij}) = \left\{
    \begin{aligned}
      &-1& &\text{if $\overleftarrow{a_{ij}} \in A$}     \\
      &0&  &\text{if $\overleftarrow{a_{ij}}, \overrightarrow{a_{ij}} \notin A$} \\
      &1&  &\text{if $\overrightarrow{a_{ij}} \in A$}&
    \end{aligned}
    \right.&
  &\text{and}&
  &\Sigma = \mathbf{O}
\end{align*}
as in the \textit{minimum entropy} case of UGs. The \textit{intermediate entropy}
case again ranges from being very close to the \textit{minimum entropy} case
(when the graph structure displays little variability) to being very close to
the \textit{maximum entropy} case (when the graph structure displays substantial
variability). The bounds on the eigenvalues of $\Sigma$ derived in Lemma
\ref{thm:dirlambda} allow a graphical representation of the variability of the
network structure, equivalent to the one illustrated in Example \ref{ex:base}
for UGs.

\section{Measures of variability}
\label{sec:variability}

Several functions have been proposed in literature as univariate measures of
spread of a multivariate distribution, usually under the assumption of multivariate
normality; for some examples see \citet{mardia} and \citet{bilodeau}. Three of
them in particular can be used as descriptive statistics for the multivariate
Bernoulli and Trinomial distributions: the \textit{generalised variance},
\begin{equation*}
  \VAR_G(\Sigma) = \det(\Sigma);
\end{equation*}
the \textit{total variance},
\begin{equation*}
  \VAR_T(\Sigma) = \tr(\Sigma);
\end{equation*}
and the squared \textit{Frobenius matrix norm} of the difference between
$\Sigma$ and a target matrix $\Psi$,
\begin{equation*}
\label{eqn:sqfrobenius}
  \VAR_F(\Sigma, \Psi) = ||| \Sigma - \Psi|||_F^2.
\end{equation*}

Both generalised variance and total variance associate high values of the
statistic to unstable network structures, and are bounded due to the properties
of the multivariate Bernoulli and Trinomial distributions. For total variance,
it is easy to show that either $\VAR_T(\Sigma) \in [0, \frac{k}{4}]$ (for the
multivariate Bernoulli) or $\VAR_T(\Sigma) \in [0, k]$ (for the multivariate
Trinomial), due to the bounds on the variances $\sigma_{ii}$ and on the eigenvalues
$\lambda_i$ derived in Sections \ref{sec:mvber} and \ref{sec:mvtri}. Generalised
variance is similarly bounded due to Hadamard's theorem on the determinant of a
non-negative definite matrix \citep{seber}: $\VAR_G(\Sigma) \in [0, (\frac{1}{4})^k]$
for the multivariate Bernoulli distribution and $\VAR_G(\Sigma) \in [0, 1]$ for
the multivariate Trinomial. They reach the respective maxima in the \textit{maximum
entropy} case and are equal to zero only in the \textit{minimum entropy} case.
Generalised variance is also strictly convex, but it is equal to zero when
$\Sigma$ is rank deficient. For this reason it may be convenient to reduce $\Sigma$
to a smaller, full rank matrix (say $\Sigma^*$) and consider $\VAR_G(\Sigma^*)$
instead of $\VAR_G(\Sigma)$; using a regularised estimator for $\Sigma$ such as
the one presented in \citet{ledoit} is also a viable option.

The behaviour of the squared Frobenius matrix norm, on the other hand, depends
on the choice of the target matrix $\Psi$. For $\Psi = \mathbf{O}$ (the covariance
matrix arising from the \textit{minimum entropy} case for both the multivariate
Bernoulli and the multivariate Trinomial), $\VAR_F(\Sigma, \Psi)$ associates high
values of the statistic to unstable network structures, like $\VAR_T(\Sigma)$ and
$\VAR_G(\Sigma)$; however, $\VAR_F(\Sigma, \mathbf{O})$ does not have a unique
maximum and none of its maxima corresponds to the \textit{maximum entropy} case,
making its interpretation unclear. A better choice seems to be a multiple of the
covariance matrix arising from the \textit{maximum entropy} case, say
$\Psi = k\Sigma_{max}$, associating high values of $\VAR_F(\Sigma, k\Sigma_{max})$
to stable network structures. For the multivariate Bernoulli, if we let
$\Psi = \frac{k}{4}I_k$, $\VAR_F(\Sigma, \frac{k}{4}I_k)$ can be rewritten as
\begin{equation*}
  \VAR_F\left(\Sigma, \frac{k}{4}I_k\right) = \sum_{i=1}^k \left( \lambda_i - \frac{k}{4}\right)^2.
\end{equation*}
It has both a unique global minimum (because it is a convex function),
\begin{equation*}
  \min_{\mathcal{L}} \VAR_F\left(\Sigma, \frac{k}{4}I_k\right) = \VAR_F\left(\frac{1}{4}I_k\right)
    = \sum_{i=1}^k \left( \frac{1}{4} - \frac{k}{4}\right)^2 = \frac{k(k-1)^2}{16},
\end{equation*}
and a unique global maximum,
\begin{equation*}
  \max_{\mathcal{L}} \VAR_F\left(\Sigma, \frac{k}{4}I_k\right) = \VAR_F(\mathbf{O})
     = \sum_{i=1}^k \left( \frac{k}{4}\right)^2 = \frac{k^3}{16},
\end{equation*}
which correspond to the \textit{maximum} and \textit{minimum entropy} covariance
matrices, respectively. Similar results can be derived for the multivariate
Trinomial distribution, using an approximate estimate for $\Sigma_{max}$
based on the results presented in Section \ref{sec:dagprop}.

All the descriptive statistics introduced in this section can be normalised as
follows:
\begin{gather}
\label{eqn:normalised}
\begin{split}
  \overline{\VAR}_T(\Sigma) = \frac{\VAR_T(\Sigma)}{\max_{\Sigma} \VAR_T(\Sigma)}, \quad 
  \overline{\VAR}_G(\Sigma) = \frac{\VAR_G(\Sigma)}{\max_{\Sigma} \VAR_G(\Sigma)}, \\
  \overline{\VAR}_F(\Sigma, k\Sigma_{max}) = \frac{\max_{\Sigma} \VAR_F(\Sigma, k\Sigma_{max}) - \VAR_F(\Sigma, k\Sigma_{max})}
    {\max_{\Sigma} \VAR_F(\Sigma, k\Sigma_{max}) - \min_{\Sigma} \VAR_F(\Sigma, k\Sigma_{max})}.
\end{split}
\end{gather}
These normalised statistics vary in the $[0,1]$ interval and associate high
values to graphs whose structures display a high variability. Since they vary
on a known and bounded scale, they are easy to interpret as absolute quantities
(i.e. goodness-of-fit statistics) as well as relative ones (i.e. proportions
of total possible variability).

They also have a clear geometric interpretation as distances in $\mathcal{L}$,
as they can all be rewritten as function of the eigenvalues $\lambda_1, \ldots,
\lambda_k$. This allows, in turn, to provide an easy interpretation of otherwise
complex properties of $\eprior$ and $\eposterior$ and to derive new results.
First of all, the measures introduced in Equation \ref{eqn:normalised} can be
used to select the best learning algorithm $\mathcal{A}$ in terms of structure
stability for a given data set $\mathcal{D}$. Different algorithms make use of
the information present in the data in different ways, under different sets of
assumptions and with varying degrees of robustness. Therefore, in practice
different algorithms learn different structures from the same data and, in
turn, result in different posterior distributions on $\mathbf{G}$. If we
rewrite Equation \ref{eqn:structlearn} to make this dependence explicit,
\begin{equation*}
  \Prob(\mathcal{G}(\mathcal{E}) \given \mathcal{D}, \mathcal{A}) \propto
    \eprior \Prob(\mathcal{D} \given \mathcal{G}(\mathcal{E}), \mathcal{A}),
\end{equation*}
and denote with $\Sigma_\mathcal{A}$ the covariance matrix of the distribution
of the edges (or the arcs) induced by $\Prob(\mathcal{G}(\mathcal{E}) \given
\mathcal{D}, \mathcal{A})$, then we can choose the optimal structure learning
algorithm $\mathcal{A}^*$ as
\begin{equation*}
  \mathcal{A}^* = \argmin_\mathcal{A} \overline{\VAR}_T(\Sigma_\mathcal{A})
\end{equation*}
or, equivalently, using $\overline{\VAR}_G(\Sigma_\mathcal{A})$ or
$\overline{\VAR}_F(\Sigma_\mathcal{A}, k\Sigma_{max})$ instead of
$\overline{\VAR}_T(\Sigma_\mathcal{A})$. Such an algorithm has the desirable
property of maximising the information gain from the data, as measured by the
distance from the non-informative prior $\eprior$ in $\mathcal{L}$. In other
words, $\mathcal{A}^*$ is the algorithm that uses the data in the most efficient
way. Furthermore, an optimal $\mathcal{A}^*$ can be identified even for data
sets without a ``golden standard'' graph structure to use for comparison; this
is not possible with the approaches commonly used in literature, which rely
on variations of Hamming distance \citep{graphs} and knowledge of such a
``golden standard'' to evaluate learning algorithms \citep[see, for example][]{mmhc}.

Similarly, it is possible to study the influence of different values of a tuning
parameter for a given structure learning algorithm (and again a given data set).
Such parameters include, for example, restrictions on the degrees of the nodes
\citep{sc} and regularisation coefficients \citep{koller}. If we denote these
tuning parameters with $\tau$, we can again choose an optimal
$\tau^*$ as
\begin{equation*}
  \tau^* = \argmin_\tau \overline{\VAR}_T(\Sigma_{\mathcal{A}(\tau)}).
\end{equation*}

Another natural application of the variability measures presented in Equation
\ref{eqn:normalised} is the study of the consistency of structure learning
algorithms. It has been proved in literature that most of structure learning
algorithms are increasingly able to identify a single, minimal graph structure
as the sample size diverges \citep[see, for example][]{hc2}. Therefore,
$\eposterior$ converges towards the \textit{minimum entropy} case and all
variability measures converge to zero. However, convergence speed has
never been analysed and compared across different learning algorithms; any
one of $\overline{\VAR}_T(\Sigma_\mathcal{A})$, $\overline{\VAR}_G(\Sigma_\mathcal{A})$
or $\overline{\VAR}_F(\Sigma_\mathcal{A}, k\Sigma_{max})$ provides a
coherent way to perform such an analysis.

Lastly, we may use the variability measures from Equation \ref{eqn:normalised}
as basis to investigate different prior distributions for real-world data
modelling and to define new ones. Relatively little attention has been paid in
literature to the choice of the prior over $\mathbf{G}$, and the uniform
\textit{maximum entropy} distribution is usually chosen for computational
reasons. Its only parameter is the \textit{imaginary sample size}, which
expresses the weight assigned to the prior distribution as the size of an
imaginary sample size supporting it \citep{heckerman}.

However, choosing a uniform prior also has some drawbacks. Firstly, \citet{steck2}
and \citet{steck1} have shown that both large and small values of the imaginary
sample size have unintuitive effects on the sparsity of a Bayesian network even
for large sample sizes. For instance, large values of the imaginary sample size
may favour the presence of an arc over its absence even when both $\eprior$ and
$\mathcal{D}$ imply the variables the arc is incident on are conditionally independent.
Secondly, a uniform prior assigns a non-null probability to all possible models.
Therefore, it often results in a very flat posterior which is not able discriminate
between networks that are well supported by the data and networks that are not
\citep{koller}. 

Following \citet{pearl}'s suggestion that ``good'' graphical models should be
sparse, sparsity-inducing priors such as the ones in \citet{buntine} and 
\citet{friedman2} should be preferred to the \textit{maximum entropy}
distribution, as should informative priors \citep{MS2008}. For example, the
prior proposed in \citet{buntine} introduces a prior probability $\beta$ to
include (independently) each arc in a Bayesian network with a given topological
ordering, which means $\overrightarrow{p_{ij}} = \beta$ and 
$\overleftarrow{p_{ij}} = 0$ for all $i < j$ in $\eprior$. Thus, $\VAR(A_{ij})
= \beta - \beta^2$, $\VAR_T(\Sigma) = k(\beta - \beta^2)$ and $\VAR_G(\Sigma)
= (\beta - \beta^2)^k$. The prior proposed in \citet{friedman}, on the other
hand, controls the number of parents of each node for a given topological
ordering. Therefore, it favours low values of $\Prob(\overrightarrow{a_{ij}},
\overleftarrow{a_{jk}})$ in $\eprior$ and again $\overleftarrow{p_{ij}} = 0$
for all $i < j$. Clearly, the amount of sparsity induced by the hyperparameters
of these priors determines the variability of both the prior and the posterior,
and can be controlled through the variability measures from Equation 
\ref{eqn:normalised}. Furthermore, these measures can provide inspiration in 
devising new priors with the desired form and amount of sparsity.

\section{Conclusions}
\label{sec:conclusion}

Bayesian inference on the structure of graphical models is challenging in most
situations due to the difficulties in defining and analysing prior and posterior
distributions over the spaces of undirected or directed acyclic graphs. The
dimension of these spaces grows super-exponentially in the number of variables
considered in the model, making even MAP analyses problematic.

In this paper, we propose an alternative approach to the analysis of graph
structures which focuses on the set of possible edges $\mathcal{E}$ of a
graphical model $\mathcal{M} = (\mathcal{G}(\mathcal{E}), \Theta)$ instead of
the possible graph structures themselves. The latter are uniquely identified by
the respective edge sets; therefore, the proposed approach integrates smoothly
with and extends both frequentist and Bayesian results present in literature.
Furthermore, this change in focus provides additional insights on the behaviour
of individual edges (which are usually the focus of inference) and reduces the
dimension of the sample space from super-exponential to quadratic in the number
of variables. 

For many inference problems the parameter space is reduced as well, and makes
complex inferential tasks feasible. As an example, we characterise several
measures of structural variability for both Bayesian and Markov networks using
the second order moments of $\eprior$ and $\eposterior$. These measures have
several possible applications and are easy to interpret from both an algebraic
and a geometric point of view.

\begin{acknowledgement}
The author would like to thank to Adriana Brogini (University of Padova) and 
David Balding (University College London) for proofreading this article and
providing many useful comments and suggestions. Furthermore, the author would
also like to thank Giovanni Andreatta and Luigi Salce (University of Padova) for
their assistance in the development of the material.
\end{acknowledgement}

\appendix

\section{Proofs}
\label{app:proofs}

\begin{proof}[Proof of Lemma \ref{thm:mvebereigen}]
  Since $\Sigma$ is a real, symmetric, non-negative definite matrix, its
  eigenvalues $\lambda_i$ are non-negative real numbers; this proves the
  lower bound in both inequalities.

  The upper bound in the first inequality holds because
  \begin{equation*}
     \sum_{i=1}^k \lambda_i = \sum_{i=1}^k \sigma_{ii} \leqslant
     \max_{\left\{\sigma_{ii}\right\}} \sum_{i=1}^k \sigma_{ii} =
     \sum_{i=1}^k \max \sigma_{ii} = \frac{k}{4},
  \end{equation*}
  as the sum of the eigenvalues is equal to the trace of $\Sigma$. This in turn
  implies
  \begin{equation*}
    \lambda_i \leqslant \sum_{i=1}^k \lambda_i \leqslant \frac{k}{4},
  \end{equation*}
  which completes the proof.
\end{proof}

\begin{proof}[Proof of Theorem \ref{thm:triber2}]
  It is easy to show that each $|T_i| = B_i$, with $p_{i(1)} + p_{i(-1)} = p^*_i$
  and $p_{i(0)} = 1 - p^*_i$. It follows that the parameter collection $\mathbf{p}$
  of $\mathbf{T}$ reduces to
  \begin{align*}
    \mathbf{p}^* &= \left\{ p_{I(T)} : I \subseteq \{1, \ldots, k\},\,
      T \in \{0, 1\}^{|I|},\, I \neq \varnothing \right\} \notag \\
      &= \left\{ p_I : I \subseteq \{1, \ldots, k\},\, I \neq \varnothing \right\}
  \end{align*}
  after the transformation. Therefore, $|\mathbf{T}| \sim Ber_k(\mathbf{p}^*)$ is
  a uniquely identified multivariate Bernoulli random variable according to the
  definition introduced at the beginning of Section \ref{sec:mvber}.
\end{proof}

\begin{proof}[Proof of Theorem \ref{thm:acyclic}]
  Let's assume by contradiction that $G^*$ is cyclic; this implies that there are
  one or more nodes $v_i \in \mathbf{V}$ such that
  \begin{equation*}
    v_i \xrightarrow{\overrightarrow{a_{ij}}} v_j \rightarrow \ldots \rightarrow v_k \xrightarrow{\overrightarrow{a_{ki}}} v_i
  \end{equation*}
  for some $v_j, v_k \in \mathbf{V}$. However, this would mean that in $G$ we would
  have
  \begin{equation*}
    v_i \xrightarrow{\overleftarrow{a_{ki}}} v_k \rightarrow \ldots \rightarrow v_j \xrightarrow{\overleftarrow{a_{ij}}} v_i
  \end{equation*}
  which is not possible since $G$ is assumed to be acyclic.
\end{proof}

\begin{proof}[Proof of Theorem \ref{thm:corollary}]
  Each possible arc can appear in the graph in only one direction at a time, so
  a directed acyclic graph with $n$ nodes can have at most ${n \choose 2} =
  \frac{1}{2}n(n - 1)$ arcs. Therefore
  \begin{equation*}
  \overrightarrow{p_{ij}} + \overleftarrow{p_{ij}} \simeq \frac{\frac{1}{4}n^2}{\frac{1}{2}n(n-1)}
    = \frac{1}{2} + \frac{1}{2(n-1)}.
  \end{equation*}
  But in the \textit{maximum entropy} case we also have that $\overrightarrow{p_{ij}}
  = \overleftarrow{p_{ij}}$, so
  \begin{align*}
    &\overrightarrow{p_{ij}} = \overleftarrow{p_{ij}} \simeq \frac{1}{4} + \frac{1}{4(n-1)}&
    &\text{and}&
    &\mathring{p_{ij}} = 1 - 2\overrightarrow{p_{ij}} \simeq \frac{1}{2} - \frac{1}{2(n-1)},
  \end{align*}
  which completes the proof.
\end{proof}

\begin{proof}[Proof of Theorem \ref{thm:hoeffding}]
  In the maximum entropy case, all arcs have the same marginal distribution function,
  \begin{equation}
  \label{eqn:distrfun}
    F_A(a_{ij}) \simeq \left\{
    \begin{aligned}
    &0                                 & & \text{in } (-\infty, -1] \\
    &\frac{1}{4} + \frac{1}{4(n - 1)}  & & \text{in } (-1, 0] \\
    &\frac{3}{4} - \frac{1}{4(n - 1)}  & & \text{in } (0, 1] \\
    &1                                 & & \text{in } (1, +\infty)
    \end{aligned}
  \right.,
  \end{equation}
  so the joint distribution of any pair of arcs $a_{ij}$ and $a_{kl}$ can be written
  as a member of the Farlie-Morgenstern-Gumbel family of distribution as
  \begin{align}
  \label{eqn:this}
    F_{A_{ij}, A_{kl}}(a_{ij}, a_{kl}) = F_A(a_{ij})F_A(a_{kl})[1 + \varepsilon(1 - F_A(a_{ij}))(1 - F_A(a_{kl}))].
  \end{align}

  Then if we apply Hoeffding's identity from Equation \ref{eqn:hoeffding} and
  replace the joint distribution function $F_{A_{ij}, A_{kl}}(a_{ij}, a_{kl})$
  with the right hand of Equation \ref{eqn:this} we have that
  \begin{align*}
    &|\COV(A_{ij}, A_{kl})| = \notag \\
    &= \left|\sum_{\{-1, 0, 1\}}\sum_{\{-1, 0, 1\}} F_{A_{ij}, A_{kl}}(a_{ij}, a_{kl}) -  F_A(a_{ij})F_A(a_{kl}) \right| \\
    &  \leqslant \sum_{\{-1, 0, 1\}}\sum_{\{-1, 0, 1\}} \left| F_{A_{ij}, A_{kl}}(a_{ij}, a_{kl}) -  F_A(a_{ij})F_A(a_{kl}) \right| \\
    &  = \sum_{\{-1, 0, 1\}}\sum_{\{-1, 0, 1\}} \left| F_A(a_{ij})F_A(a_{kl})[1 + \right. \\
    & \qquad\qquad\qquad\qquad\qquad \left. + \varepsilon(1 - F_A(a_{ij}))(1 - F_A(a_{kl}))]  -  F_A(a_{ij})F_A(a_{kl}) \right| \\
    &  = \sum_{\{-1, 0\}}\sum_{\{-1, 0\}} (1 - F_A(a_{ij}))(1 - F_A(a_{kl})). 
  \end{align*}

  We can now compute the bounds for $|\COV(a_{ij}, a_{kl})|$ and $|\COR(a_{ij}, a_{kl})|$
  using only the marginal distribution function $F_A$ from Equation \ref{eqn:distrfun}
  and the variance from Equation \ref{eqn:approxvar}, thus obtaining the
  expressions in Equation \ref{eqn:boundcov} and Equation \ref{eqn:boundcor}.
\end{proof}

\section{Moments and parameters of the multivariate Trinomial distribution in
  the maximum entropy case}
\label{app:numbers}

Below are reported the exact values of the parameters of the marginal Trinomial
distributions and of the first and second order moments of the multivariate
Trinomial distribution in the maximum entropy case. All these quantities have
been computed by a complete enumeration of the directed acyclic graphs of a
given size ($3$, $4$, $5$, $6$ and $7$).

\subsection{Moments for the 3-dimensional distribution}

\begin{align*}
  &A_{ij} = \left\{
    \begin{aligned}
      &-1& &\text{with probability $0.32$} \\
      &0&  &\text{with probability $0.36$} \\
      &1&  &\text{with probability $0.32$}&
    \end{aligned}
    \right.&
  & \begin{aligned}
      &\E(A_{ij}) = 0 \\
      &\VAR(A_{ij}) = 0.64 \\
      &|\COV(A_{ij}, A_{kl})| = 0.08
    \end{aligned}
\end{align*}

\subsection{Moments for the 4-dimensional distribution}

\begin{align*}
  &A_{ij} = \left\{
    \begin{aligned}
      &-1& &\text{with probability $0.309392$} \\
      &0&  &\text{with probability $0.381215$} \\
      &1&  &\text{with probability $0.309392$}&
    \end{aligned}
    \right.&
  & \begin{aligned}
      &\E(A_{ij}) = 0 \\
      &\VAR(A_{ij}) = 0.618784 
    \end{aligned}
\end{align*}

\begin{equation*}
  |\COV(A_{ij}, A_{kl})| = \left\{
    \begin{aligned}
      &0&  &\text{if $i \neq j \neq k \neq l $} \\
      &0.081031&  &\text{otherwise}&
    \end{aligned}
  \right.
\end{equation*}

\subsection{Moments for the 5-dimensional distribution}

\begin{align*}
  &A_{ij} = \left\{
    \begin{aligned}
      &-1& &\text{with probability $0.301082$} \\
      &0&  &\text{with probability $0.397834$} \\
      &1&  &\text{with probability $0.301082$}&
    \end{aligned}
    \right.&
  & \begin{aligned}
      &\E(A_{ij}) = 0 \\
      &\VAR(A_{ij}) = 0.602165
    \end{aligned}
\end{align*}

\begin{equation*}
  |\COV(A_{ij}, A_{kl})| = \left\{
    \begin{aligned}
      &0&  &\text{if $i \neq j \neq k \neq l $} \\
      &0.081691&  &\text{otherwise}&
    \end{aligned}
  \right.
\end{equation*}

\subsection{Moments for the 6-dimensional distribution}

\begin{align*}
  &A_{ij} = \left\{
    \begin{aligned}
      &-1& &\text{with probability $0.294562$} \\
      &0&  &\text{with probability $0.410875$} \\
      &1&  &\text{with probability $0.294562$}&
    \end{aligned}
    \right.&
  & \begin{aligned}
      &\E(A_{ij}) = 0 \\
      &\VAR(A_{ij}) = 0.589124
    \end{aligned}
\end{align*}

\begin{equation*}
  |\COV(A_{ij}, A_{kl})| = \left\{
    \begin{aligned}
      &0&  &\text{if $i \neq j \neq k \neq l $} \\
      &0.082121&  &\text{otherwise}&
    \end{aligned}
  \right.
\end{equation*}

\subsection{Moments for the 7-dimensional distribution}

\begin{align*}
  &A_{ij} = \left\{
    \begin{aligned}
      &-1& &\text{with probability $0.289390$} \\
      &0&  &\text{with probability $0.421220$} \\
      &1&  &\text{with probability $0.289390$}&
    \end{aligned}
    \right.&
  & \begin{aligned}
      &\E(A_{ij}) = 0 \\
      &\VAR(A_{ij}) = 0.578780 
    \end{aligned}
\end{align*}

\begin{equation*}
  |\COV(A_{ij}, A_{kl})| = \left\{
    \begin{aligned}
      &0&  &\text{if $i \neq j \neq k \neq l $} \\
      &0.82410&  &\text{otherwise}&
    \end{aligned}
  \right.
\end{equation*}

\bibliographystyle{ba}

\end{document}